\documentclass{article}

\usepackage{lineno,hyperref}
\usepackage{amssymb}
\usepackage{amsmath,amsfonts,tikz,graphicx,subfig,amsthm}
\usepackage{graphicx}

\newsavebox{\bigpicture}
\newcommand{\adapttobigpicture}{%
  \vrule height\ht\bigpicture depth\dp\bigpicture width0pt
}

\newtheorem{mydef-theorem}{Theorem}
\newtheorem{mydef-observation}{Observation}
\newtheorem{mydef-propostion}{Propostion}
\newtheorem{mydef-corollary}{Corollary}
\date{}

\begin{document}

\title{Structural and Spectral properties of Corona Graphs}

\author{Rohan Sharma \, \thanks{Centre for System Science, Indian Institute of Technology Jodhpur, Jodhpur, India, Email id: rohan.sharma@iitj.ac.in} Bibhas Adhikari \, \thanks{ Department of Mathematics, Indian Institute of Technology Kharagpur, Kharagpur, India, Email id: bibhas@maths.iitkgp.ernet.in} and \, Abhishek Mishra\thanks{Department of Computer Science and Information Systems, BITS-Pilani, Pilani, India}}

\maketitle

\begin{abstract}
Product graphs have been gainfully used in literature to generate mathematical models of complex networks which inherit properties of real networks. Realizing the duplication phenomena imbibed in the definition of corona product of two graphs, we define Corona graphs. Given a small simple connected graph which we call seed graph, Corona graphs are defined by taking corona product of a seed graph iteratively. We show that the cumulative degree distribution of Corona graphs decay exponentially when the seed graph is regular and cumulative betweenness distribution follows power law when seed graph is a clique. We determine explicit formulae of eigenvalues, Laplacian eigenvalues and signless Laplacian eigenvalues of Corona graphs when the seed graph is regular. Computable expressions of eigenvalues and signless Laplacian eigenvalues of Corona graphs are also obtained when the seed graph is a star graph.
\end{abstract}

Keywords: Corona product of graphs, Laplacian spectra, signless Laplacian spectra , complex networks , betweenness distribution , cumulative degree distribution
\section{Introduction}
Network modelling using product graphs is an interesting method to generate complex networks which possibly can capture the properties of  real world networks, for example, Kronecker graphs (see \cite{leskovec2010kronecker}). In \cite{leskovec2010kronecker}, Leskovec et.al. have developed both deterministic and stochastic models by using the Kronecker product of graphs to generate complex networks which inherit properties of real world networks. In \cite{parsonage2011generalized}, the authors had presented generalized graph products methodology to generate different types of complex networks. In this paper, we propose a model for complex networks by using the concept of corona product of graphs.

\begin{figure}
\centering
\sbox{\bigpicture}{%

\begin{tikzpicture}
  [scale=0.5,auto=left,every node/.style={circle,fill=black!100}]
  \node (n1) at (0,0) {};
  \node (n2) at (1,-1) {};
  \node (n3) at (-1,-1) {};
  
  \foreach \from/\to in {n1/n2,n1/n3}
    \draw (\from) -- (\to);
\end{tikzpicture}
}
\subfloat[$P_{3}$]{\usebox{\bigpicture}\label{1a}}\qquad
  \subfloat[$G^{(1)}$]{%
  \adapttobigpicture
\begin{tikzpicture}
  [scale=0.5,auto=left,every node/.style={circle,fill=black!100}]
  \node (n1) at (0,0) {};
  \node (n2) at (1,-1) {};
  \node (n3) at (-1,-1) {};
  \node (n4) at (0,2) {};
  \node (n5) at (-1,1) {};
  \node (n6) at (1,1) {};
  \node (n7) at	(-2,-0.5) {};
  \node (n8) at (-3,-2) {};
  \node (n9) at (-1,-2) {};
  \node (n10) at (2,-0.5) {};
  \node (n11) at (1,-2) {};
  \node (n12) at (3,-2) {};
  \foreach \from/\to in {n1/n2,n1/n3,n1/n4,n1/n5,n1/n6,n4/n5,n4/n6,n2/n10,n2/n11,n2/n12,n10/n12,n11/n12, n3/n7,n3/n9,n3/n8,n7/n8,n8/n9}
    \draw (\from) -- (\to);
\end{tikzpicture}
\label{1b}
}
\subfloat[$G^{(2)}$]
{%
  \adapttobigpicture
  \begin{tikzpicture}
    [scale=0.5,auto=left,every node/.style={circle,fill=black!100}]
    \node (n1) at (0,0) {};
    \node (n2) at (1,-1) {};
    \node (n3) at (-1,-1) {};
    \node (n4) at (0,2) {};
    \node (n5) at (-1,1) {};
    \node (n6) at (1,1) {};
    \node (n7) at	(-2,-0.5) {};
    \node (n8) at (-3,-2) {};
    \node (n9) at (-1.5,-2) {};
    \node (n10) at (2,-0.5) {};
    \node (n11) at (1.5,-2) {};
    \node (n12) at (3,-2) {};
    \node (n13) at (-2,3) {};
    \node (n14) at (-2.75,2) {};
    \node (n15) at (-0.8,2){};
    \node (n16) at (0,5) {};
    \node (n17) at (-1,4) {};
    \node (n18) at (1,4) {};
    \node (n19) at (2.5,2) {};
    \node (n20) at (0.7,2) {};
    \node (n21) at (1.75,3) {};
    \node (n22) at (4,0.75) {};
    \node (n23) at (4,2) {};
    \node (n24) at (5,1.5) {};
    \node (n25) at (4,0) {};
    \node (n26) at (5,-1) {};
    \node (n27) at (4,-1.5) {};
    \node (n28) at (5,-2.25) {};
    \node (n29) at (6,-3) {};
    \node (n30) at (5,-3.5) {};
    \node (n31) at (3.5,-3.5) {};
    \node (n32) at (3.75,-5) {};
    \node (n33) at (2.5,-4.5) {};
    \node (n34) at (1.15,-6) {};
    \node (n35) at (1.75,-5) {};
    \node (n36) at (0.5,-5) {};
    \node (n37) at (-4,0) {};
    \node (n38) at (-5,-1) {};
    \node (n39) at (-4,-1.5) {};
    \node (n40) at (-5,-2.05) {};
    \node (n41) at (-6,-3) {};
    \node (n42) at (-5,-3.5) {};
    \node (n43) at (-4.05,-3.5) {};
    \node (n44) at (-4.25,-5.05) {};
    \node (n45) at (-2.75,-4.5) {};
    \node (n46) at (-1,-3.5) {};
    \node (n47) at (0,-4) {};
    \node (n48) at (0.25,-2.5) {};                
    \foreach \from/\to in {n1/n2,n1/n3,n1/n4,n1/n5,n1/n6,n4/n5,n4/n6,n2/n10,n2/n11,n2/n12,n10/n12,n11/n12,n3/n7,n3/n9,n3/n8,n7/n8,n8/n9,n5/n13,n5/n14,n5/n15,n13/n14,n15/n13,n17/n16,n16/n18,n4/n16,n4/n17,n4/n18,n21/n19,n21/n20,n6/n19,n6/n20,n6/n21,n1/n22,n1/n23,n1/n24,n22/n24,n24/n23,n25/n26,n26/n27,n10/n25,n10/n26,n10/n27,n28/n29,n30/n29,n28/n12,n12/n29,n12/n30,n31/n32,n32/n33,n11/n31,n11/n32,n11/n33,n2/n35,n2/n36,n2/n34,n35/n34,n36/n34,n37/n38,n39/n38,n7/n38,n37/n7,n7/n39,n40/n8,n41/n8,n42/n8,n40/n41,n41/n42,n9/n43,n9/n44,n9/n45,n43/n44,n44/n45,n3/n46,n3/n47,n3/n48,n46/n47,n47/n48}
      \draw (\from) -- (\to);
  \end{tikzpicture}
\label{1c}
  }
\caption{\label{1}Examples of the Corona graphs: (a) A seed graph $P_{3}$ (b) Corona product of $P_3=G$ with itself and hence generating $G^{(1)}$ (c) The graph of another successive corona product of $P_3$ with previously generated graph $G^{(1)}$ resulting in $G^{(2)}$.}
\end{figure}
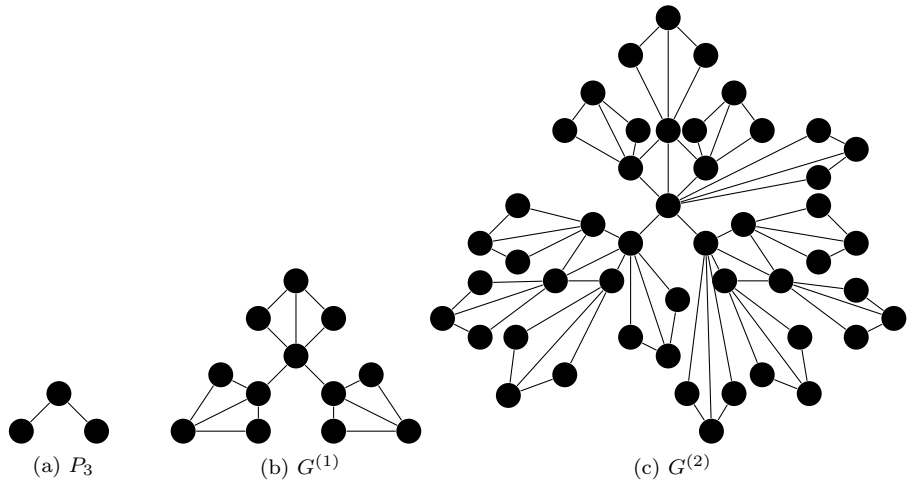

Corona product of two graphs, say $G$ and $H$, was introduced by Frucht and Harary in 1970 \cite{frucht1970corona} is a graph constructed by taking $n$ instances of $H$ and each such $H$ gets connected to each node of $G,$ where $n$ is the number of nodes of $G.$ Starting with a connected simple graph $G,$ we define Corona graphs which are obtained by taking Corona product of $G$ with itself iteratively. In this case, $G$ is called the seed graph for the Corona graphs. For instance, the Corona graphs generated with $G=P_3$ is shown in Fig.~\ref{1}.

We mention that Corona graphs can be used as models for investigating the duplication mechanism of an individual gene as explained for the formation of new proteins in \cite{ispolatov2005duplication} and the references therein. A similar phenomenon is being observed in Corona graphs but instead of a single node, a unit of seed graph is being duplicated and connected to every node of the existing network. Hence, Corona graphs can reveal more insights on duplication phenomena and with the help of proposed graph spectra, the properties of gene duplication and other real world complex networks can be investigated more deeply \cite{parsonage2011generalized}.

In this paper, we investigate various structural properties of Corona graphs including average degree, sparsity, degree distribution, diameter and cumulative degree distribution. We show that \begin{enumerate}
\item the cumulative degree distribution of a Corona graph decays exponentially when the chosen seed graph is regular. \item diameter of the Corona graphs is increased by $2$ in each iteration of the Corona product for any seed graph. \item cumulative betweenness distribution of a Corona graph follows power law when  the seed graph is a clique.
\end{enumerate}

Further, we study the spectra, Laplacian spectra and signless Laplacian spectra of Corona graphs obtained by specific seed graphs. We determine
\begin{enumerate}
 \item computable formulae for eigenvalues and signless Laplacian eigenvalues of Corona graphs generated by a seed graph which is regular or a star graph. \item computable expressions of Laplacian eigenvalues of Corona graphs generated by any seed graph.
\end{enumerate}

We organize the paper as follows. In Sec.~\ref{sec:CG}, we define the Corona graphs and investigate structural properties of Corona graphs. In Sec.~\ref{sec:Spec}, we derive the spectrum of Corona graphs generated by a regular graph and a star graph. Further we determine formula of Laplacian and signless Laplacian eigenvalues of Corona graphs generated by a simple connected graph and finally, we conclude in Sec.~\ref{sec:Concl}.
\section{Corona Graphs}\label{sec:CG}
Let $G=(V,E)$ be a graph having the set of nodes $V=\{v_1,...,v_n\}$ and E the set of edges. The adjacency matrix $A(G)=[a_{v_i v_j}]$ of $G$ having dimension $|V|\times |V|$ is defined by $a_{v_i,v_j}=1$ if $(v_i,v_j)\in E$ otherwise $a_{v_i,v_j}=0$ when columns and rows are labelled by nodes of $G$. Laplacian matrix and signless Laplacian matrix associated with $G$ are given by
\begin{equation}\label{eqn:1}
L(G)=D(G)-A(G)
\end{equation}
and
\begin{equation}\label{eqn:2}
Q(G)=D(G)+A(G)
\end{equation}
respectively, where $D(G)=diag\{d_1,d_2,\ldots,d_n\}$, $d_i=\sum_{i\neq j}a_{ij}$.
\begin{figure}
\centering

  \subfloat[$G^{(1)}$]
  {%
    \adapttobigpicture
   \includegraphics[width=0.22\textwidth]{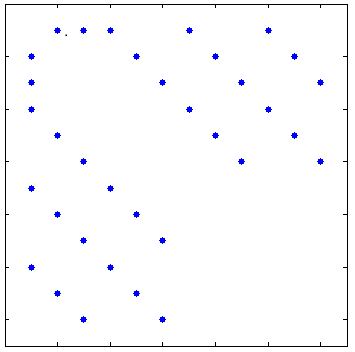} 
   \label{2a}
  }
  \subfloat[$G^{(3)}$]
  {%
    \adapttobigpicture
   \includegraphics[width=0.22\textwidth]{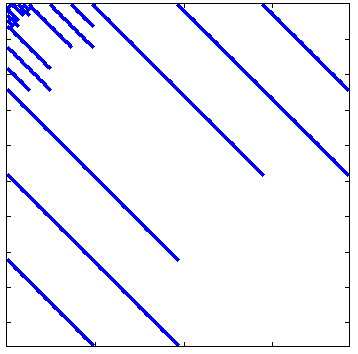} 
      \label{2b}
  }
\caption{Pattern of the adjacency matrices corresponding to (a) $G^{(1)}$ (b) $G^{(3)}$ for Fig.~\ref{1a}. Dots represent $a_{v_i,v_j}=1$ and white spaces represents $a_{v_i,v_j}=0$.} 
\end{figure}
\par
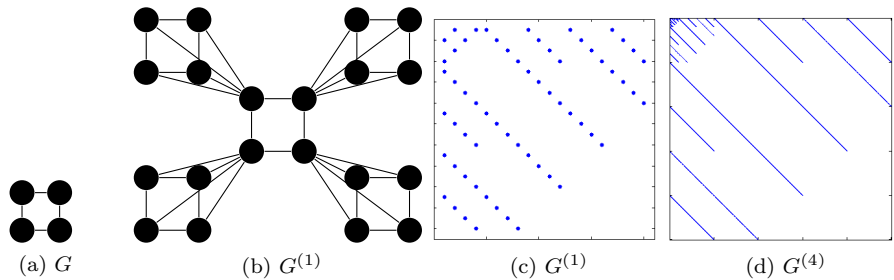
\begin{figure}
\centering
\sbox{\bigpicture}{%

\begin{tikzpicture}
  [scale=0.5,auto=left,every node/.style={circle,fill=black!100}]
  \node (n1) at (-0.5,-0.5) {};
  \node (n2) at (0.5,-0.5) {};
  \node (n3) at (0.5,0.5) {};
  \node (n4) at (-0.5,0.5) {};
  \foreach \from/\to in {n1/n2,n2/n3,n3/n4,n4/n1}
    \draw (\from) -- (\to);
\end{tikzpicture}
\label{3a}
}
\subfloat[$G$]{\usebox{\bigpicture}\label{$G$}}\qquad
  \subfloat[$G^{(1)}$]{%
  \adapttobigpicture
\begin{tikzpicture}
  [scale=0.7,auto=left,every node/.style={circle,fill=black!100}]
  \node (n1) at (-0.5,-0.5) {};
  \node (n2) at (0.5,-0.5) {};
  \node (n3) at (0.5,0.5) {};
  \node (n4) at (-0.5,0.5) {};
  \node (n5) at (1.5,1) {};
  \node (n6) at (2.5,1) {};
  \node (n7) at (2.5,2) {};
  \node (n8) at (1.5,2) {};
  \node (n9) at (1.5,-1) {};
  \node (n10) at (2.5,-1) {};
  \node (n11) at (2.5,-2) {};
  \node (n12) at (1.5,-2) {};
  \node (n13) at (-1.5,-1) {};
  \node (n14) at (-2.5,-1) {};
  \node (n15) at (-2.5,-2) {};
  \node (n16) at (-1.5,-2) {};
  \node (n17) at (-1.5,1) {};
  \node (n18) at (-2.5,1) {};
  \node (n19) at (-2.5,2) {};
  \node (n20) at (-1.5,2) {};
  \foreach \from/\to in {n1/n2,n2/n3,n3/n4,n4/n1,n5/n6,n6/n7,n7/n8,n8/n5,n9/n10,n10/n11,n11/n12,n12/n9,n13/n14,n14/n15,n15/n16,n16/n13,n17/n18,n18/n19,n19/n20,n20/n17,n1/n13,n14/n1,n16/n1,n15/n1,n9/n2,n10/n2,n11/n2,n12/n2,n5/n3,n6/n3,n7/n3,n8/n3,n4/n17,n4/n18,n4/n19,n4/n20}
    \draw (\from) -- (\to);
\end{tikzpicture}
\label{3b}
}
\subfloat[$G^{(1)}$]
  {%
    \adapttobigpicture
   \includegraphics[width=0.25\textwidth]{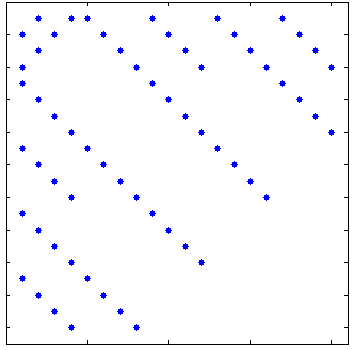} 
   \label{3c}
  }
  \subfloat[$G^{(4)}$]
    {%
      \adapttobigpicture
     \includegraphics[width=0.25\textwidth]{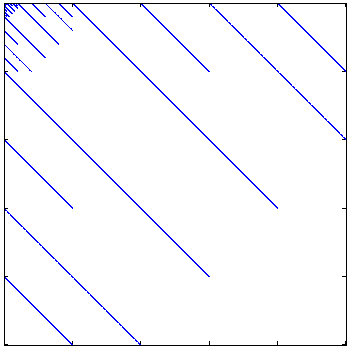} 
     \label{3d}
    }
  
\caption{Examples of the Corona graphs: (a) $2$-regular seed graph $G$ (b) Corona graph for $G^{(1)}$. Pattern of the adjacency matrices in (c) $G^{(1)}$ (d) $G^{(4)}$.}
\end{figure}

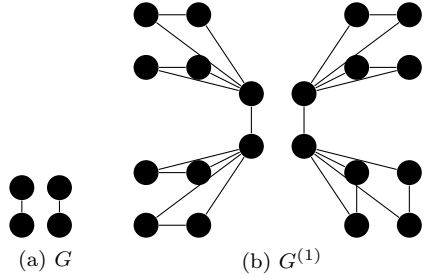
\begin{figure}
\centering
\sbox{\bigpicture} {%
\begin{tikzpicture}
  [scale=0.5,auto=left,every node/.style={circle,fill=black!100}]
  \node (n1) at (-0.5,-0.5) {};
  \node (n2) at (0.5,-0.5) {};
  \node (n3) at (0.5,0.5) {};
  \node (n4) at (-0.5,0.5) {};
  \foreach \from/\to in {n2/n3,n4/n1}
    \draw (\from) -- (\to);
\end{tikzpicture}
\label{4a}
}
\subfloat[$G$]{\usebox{\bigpicture}\label{$G$}}\qquad
  \subfloat[$G^{(1)}$] {%
  \adapttobigpicture

\begin{tikzpicture}
  [scale=0.7,auto=left,every node/.style={circle,fill=black!100}]
  \node (n1) at (-0.5,-0.5) {};
  \node (n2) at (0.5,-0.5) {};
  \node (n3) at (0.5,0.5) {};
  \node (n4) at (-0.5,0.5) {};
  \node (n5) at (1.5,1) {};
  \node (n6) at (2.5,1) {};
  \node (n7) at (2.5,2) {};
  \node (n8) at (1.5,2) {};
  \node (n9) at (1.5,-1) {};
  \node (n10) at (2.5,-1) {};
  \node (n11) at (2.5,-2) {};
  \node (n12) at (1.5,-2) {};
  \node (n13) at (-1.5,-1) {};
  \node (n14) at (-2.5,-1) {};
  \node (n15) at (-2.5,-2) {};
  \node (n16) at (-1.5,-2) {};
  \node (n17) at (-1.5,1) {};
  \node (n18) at (-2.5,1) {};
  \node (n19) at (-2.5,2) {};
  \node (n20) at (-1.5,2) {};
  \foreach \from/\to in {n2/n3,n4/n1,n5/n6,n7/n8,n10/n11,n12/n9,n13/n14,n15/n16,n17/n18,n19/n20,n1/n13,n14/n1,n16/n1,n15/n1,n9/n2,n10/n2,n11/n2,n12/n2,n5/n3,n6/n3,n7/n3,n8/n3,n4/n17,n4/n18,n4/n19,n4/n20}
    \draw (\from) -- (\to);

\end{tikzpicture}
\label{4b}
}
\caption{A disconnected graph $G$ of 4 nodes having two $P_2$ graphs as its components, along with its $1^{st}$ corona product $G^{(1)}$ resulting in a disconnected Corona graph.}
\label{4}
\end{figure}

 The corona product of the two graphs $G_1=(V(G_1),E(G_1))$ and $G_2=(V(G_2),E(G_2))$, denoted by $G_1\circ G_2$ is obtained by taking an instance of $G_1$ and $|V_{G_1}|$ instances of $G_2$ and hence connecting the $i^{th}$ node of $G_1$ to every node in the $i^{th}$ instance of $G_2$ for each $i$ \cite{barik2007spectrum}. We extend this definition to define Corona graphs. Let $G^{(0)}= G$. Given a seed graph $G,$ the Corona graphs generated by $G$ are defined by
\begin{equation}\label{eqn:3}
G^{(m+1)}=G^{(m)}\circ G
\end{equation}
where $m(\ge 0)$ is a large natural number. For instance, the Corona graphs generated by $P_3$ are shown in Fig.\ref{1b} and Fig.\ref{1c} along with the pattern of their adjacency matrices in Fig.\ref{2a} and Fig.\ref{2b} respectively. Similarly, for a $2$-regular graph with $4$ nodes, the Corona graph $G^{(1)}$ corresponding to $G$ is shown in Fig.\ref{3b} and the pattern of their adjacency matrices of $G^{(1)}$ and $G^{(4)}$ in Fig.\ref{3c} and Fig.\ref{3d} respectively. 
\noindent The following are some observations associated with Corona graphs generated by a seed graph of order $n$.
\begin{enumerate}
\item 
The number of nodes in $G^{(m)}$ is
\begin{equation}\label{eqn:4}
|V^{(m)}|=n(n+1)^{m}.
\end{equation}
\item If $|E|$ is the number of edges in the seed graph $G$, then the number of edges in $G^{(m)}$ is 
\begin{equation}\label{eqn:5}
|E^{(m)}|=(|E| + (|E| + n)((n+1)^{m}-1)).
\end{equation}
\item The number of nodes added in $i^{th}(i\le m)$ step of the formation of $G^{(m)}$ is $n^{2}(n+1)^{i-1}.$
\item Connectivity of $G^{(m)}$: Since $G$ is connected, Corona graphs generated by $G$ are connected graphs. Evidently, if seed graph $G$ is disconnected, it generates disconnected Corona graphs as shown in Fig.\ref{4}. This feature is not observed in the case of Kronecker graphs given by Leskovec et al. in \cite{leskovec2010kronecker} where the authors had taken a multi-graph as a seed graph to ensure the connectivity of $G^{(m)}$.
\item Degree sequence of corona graphs: Assume that degree sequence of $G^{(0)}= G$ is given by \{$d_{i_1}^{(0)}$, $d_{i_2}^{(0)}$,\ldots,$d_{i_n}^{(0)}$\} where $d_{i_l}^{(j)}$ represents the degree of node $i_l$ at $j^{th}$ corona product, and $x$ is the total distinct degrees. The degree sequence of $G^{(m)}$ is given by $\{D_{i_j}^{(1)}$, $(D_{i_j}^{(2)},$ $\ldots,n$ times$),\ldots,$ $(D_{i_j}^{(x)},$\ldots,$n(n+1)^{m-1}$ times$)\}$,
where $D_{i_j}^{(1)}=d_{i_j}+mn$, $D_{i_j}^{(2)}=d_{i_j}+(m-1)n$,\ldots, $D_{i_j}^{(x)}=d_{i_j}+1$. 
 \end{enumerate}

\begin{mydef-propostion}
Let $G^{(m)}=G^{(m-1)}\circ G$ be the Corona graph generated by the seed graph $G$. Let $|E|$ and $|V|=n$ be the number of edges and nodes in the seed graph $G$. Let $|E^{(m)}|$ and $|V^{(m)}|$ be the number of edges and nodes in $G^{(m)}$. Then, average degree($\langle k\rangle$) of a node in $G^{(m)}$ is
\begin{equation*}
\langle k\rangle\approx 2\biggl(1+\frac{|E|}{n}\biggr).
\end{equation*}
\end{mydef-propostion}
\begin{proof}
By the definition of average degree of a node, we have
\begin{eqnarray*}
\langle k\rangle &=&\frac{2|E^{(m)}|}{|V^{(m)}|}\\
				 &=&\frac{2((|E| + (|E| + n)((n+1)^{m}-1)))}{n(n+1)^{m}}\\
				 &=&\frac{2(|E| + n)}{n}-\frac{2}{(n+1)^{m}}\\ 
				 \mbox{If $m\gg 1$, then $\frac{2}{(n+1)^{m}}\rightarrow 0$. Hence,}\\
\langle k\rangle &\approx & 2\biggl(1+\frac{|E|}{n}\biggr).
\end{eqnarray*}
\end{proof}
It follows from above proposition that if the seed graph is a tree such that $|E|=(n-1)$, then $\langle k\rangle=2\bigl(1+\frac{n-1}{n}\bigr)=2(2-\frac{1}{n})$. If the seed graph is a clique ($K_n$, having $|E|=\binom{n}{2}$ edges and $n\ge 3$), then $\langle k\rangle\approx (n+1).$

Density of a network is an important concept to determine the sparsity of the network. The density of an undirected network $G=(V,E)$ is defined by \cite{laurienti2011universal}
\begin{equation}\label{eqn:6}
d=\frac{|E|}{\binom{n}{2}}
\end{equation}
where $|E|$ and $|V|=n$ represent the number of edges and nodes respectively in $G$. It is evident that $0<d\le 1.$ If $d\ll 1$, then the network is called sparse. It is noted in \cite{i2003optimization} that many real world networks have $d\in[10^{-5},10^{-1}]$. In \cite{zhou2006hierarchical}, the authors remarked by citing the example of neural networks that sparsity saves the energy of a network without effecting its functionality. The network sparsity is also observed in other biological networks like metabolic network (\cite{wagner2001small} and the references therein). 

\begin{mydef-propostion}
Let $|E|$ and $n$ be the number of edges and nodes respectively in seed graph $G$. Let $|E^{(m)}|$ and $|V^{(m)}|$ be the number of edges and nodes respectively in $G^{(m)}=G^{(m-1)}\circ G$.  Then, for $m\gg 1$, the Corona graphs $G^{(m)}$ are sparse as $d\ll 1$.
\end{mydef-propostion}
\begin{proof} Using equation (~\ref{eqn:6}), the density of $G^{(m)}$ is given by
\begin{eqnarray*}
d &=&\frac{|E^{(m)}|}{\binom{|V^{(m)}|}{2}}\\
  &=&\frac{2(|E|+n)}{n(n(n+1)^m-1)}-\frac{n}{\binom{n(n+1)^m}{2}}.  
\end{eqnarray*}
If $m\gg 1$, then $\frac{n}{\binom{n(n+1)^m}{2}}\rightarrow 0$ and also $\frac{2(|E|+n)}{n(n(n+1)^m-1)}\ll 1.$ Hence, the result follows.
\end{proof}

\subsection{Degree Distribution}\label{sec:Others}
Now, we consider the degree distribution of Corona graphs. In a network $G$, let $P(k)$ denote the fraction of vertices having degree $k$ and hence the degree distribution is the probability distribution of the different degrees in the whole network \cite{newman2003structure},\cite{albert2002statistical}. The total instances for a particular degree $k$ in $G^{(m)}$ is given by
\begin{equation}\label{eqn:7}
N_k=\sum_{j=1}^{n}\biggl(\delta_{k,D_{i_j}^{(1)}}+n\delta_{k,D_{i_j}^{(2)}}+n(n+1)\delta_{k,D_{i_j}^{(3)}}+\ldots+n(n+1)^{(m-1)}\delta_{k,D_{i_j}^{(x)}}\biggr)
\end{equation}
where $\delta_{k,D_{i_j}^{(b)}}$ is the Kronecker delta function. Therefore, the degree distribution $P(k)$  of Corona graphs is given by
\begin{equation}\label{eqn:8}
P(k) =\frac{\sum_{j=1}^{n}\biggl(\delta_{k,D_{i_j}^{(1)}}+n\delta_{k,D_{i_j}^{(2)}}+n(n+1)\delta_{k,D_{i_j}^{(3)}}+\ldots+n(n+1)^{(m-1)}\delta_{k,D_{i_j}^{(x)}}\biggr)}{n(n+1)^m}
\end{equation}
\par
\begin{figure}
\centering

  \subfloat[For $K_3$ with $G^{(6)}$ and $G^{(7)}$]
  {%
    \adapttobigpicture
   \includegraphics[width=0.45\textwidth]{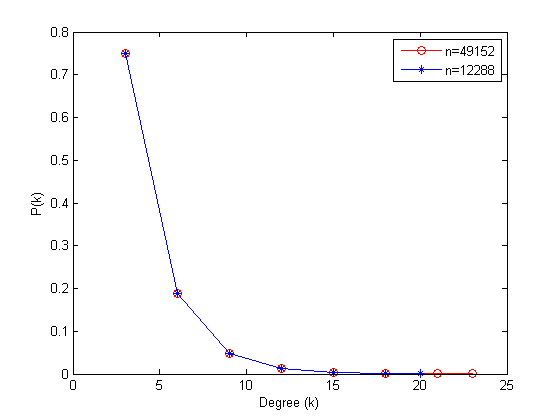} 
   \label{5a}
  }
  \subfloat[For $P_3$ with $G^{(6)}$ and $G^{(7)}$]
  {%
    \adapttobigpicture
   \includegraphics[width=0.45\textwidth]{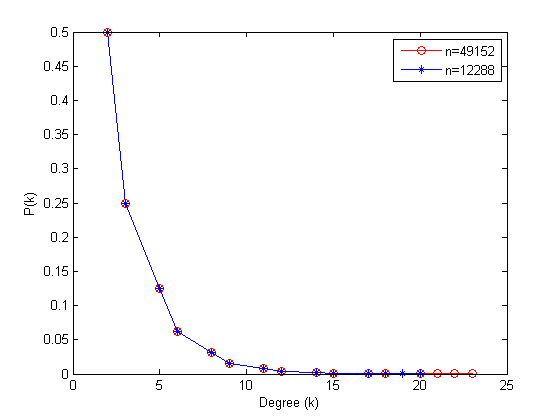} 
      \label{5b}
  }
\caption{Degree Distribution for seed graphs (a) $K_3$ with $G^{(7)}$ and $G^{(6)}$ having $49152$ and $12288$ nodes respectively. (b) $P_3$ with $G^{(7)}$ and $G^{(6)}$ having $49152$ and $12288$ nodes respectively.}
   \label{5}
\end{figure}

 The degree distribution for $K_3$ and $P_3$ are shown in Fig.\ref{5}. A crucial observation from the degree distribution is that there is huge multiplicities of degrees in the Corona graphs and this can also be confirmed from the Fig.\ref{5a}, Fig.\ref{5b}. It is similar to the observation as seen for Kronecker graphs \cite{leskovec2010kronecker}. The figure also shows that the degree distribution follows similar type of curves for different instances of Corona graph i.e. for $G^{(i)}$ for each $i\in [1,m]$. The figure also shows the fat tailed degree distribution in both sub-figures of Fig.\ref{5}. The fat tailed distributions are found abundantly in real world networks like data traffic on internet, return on financial markets etc. as mentioned in \cite{misiewicz2011fat} and the references therein \cite{crovella1999estimating},\cite{rachev2003handbook}. It should be noted that the degree sequence of a Corona graph comprises of discrete values, so the degree distribution of a Corona graph could be properly investigated with cumulative degree distribution.
\begin{equation*}
P_c(k)=\sum\limits_{k^{'}\ge k}P(k^{'})
\end{equation*}
which is the probability that the degree is greater than or equal to $k$. Here, $P(k)$ and $P_c(k)$ are the discrete and continuous degree distribution respectively. Thus, we have the following result.
\begin{mydef-theorem}\label{thm:expo}
Consider the Corona graph $G^{(t)}=G^{(t-1)}\circ G$ generated by the seed graph $G$ which is $r$-regular. Then, the cumulative degree distribution of Corona graphs decay exponentially. 
\end{mydef-theorem}
\begin{proof}
The generating methodology of Corona graphs $G^{(t)}$ for $t\ge 1$ implies
\begin{equation}\label{eqn:9}
k_i(t)=k_i(t-1)+n
\end{equation}
where $k_i(t)$ is the degree of the node $i$ in $G^{(t)}.$ Since the degree of any node is $r+1$ while it gets added during formation of $G^{(t)}$, it implies that
\begin{equation}\label{eqn:10}
k_i(t)=r+1+n(t-t_i).
\end{equation}
Now, cumulative degree distribution can be obtained using equation (~\ref{eqn:10}) as follows
\begin{equation*}
P_c(k)=\sum\limits_{k^{'}\ge k}P(k^{'}) = P(t^{'}\le \tau = t+\frac{r+1-k}{n}).
\end{equation*}
Hence,
\begin{eqnarray}
P_c(k)&=& \sum\limits_{t^{'}=0}^{\tau}\frac{|V^{(t^{'})}|}{|V^{(t)}|}\nonumber\\
	  &=& \frac{n}{|V^{(t)}|}+\sum\limits_{t^{'}=1}^{\tau}\frac{|V^{(t^{'})}|}{|V^{(t)}|}\nonumber\\
	  &=& \frac{1}{(n+1)^{t}}+\frac{((n+1)^{\tau}-1)}{(n+1)^t}\nonumber\\
	  &=& (n+1)^{\tau-t}\nonumber\\
	  \mbox{Since $\tau=t+\frac{r+1-k}{n}$, hence}\nonumber\\
P_c(k)&=& (n+1)^{\frac{r+1-k}{n}}.\label{eqn:11}
\end{eqnarray}
Since $k \geq r+1,$ the result follows.
\end{proof}

\begin{figure}
\centering
     \includegraphics[width=0.5\textwidth]{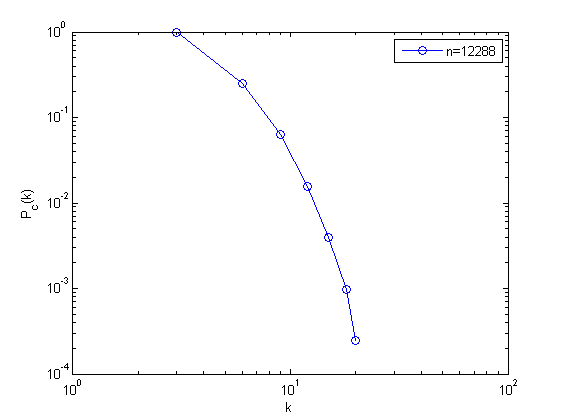} 
\captionsetup{justification=raggedright,singlelinecheck=false}
\caption{\label{fig:6}Cumulative degree distribution of $G^{(6)}$ generated with the seed graph $G^{(0)}=K_3$ where $k$ is degree of nodes and $P_c(k)$ is cumulative probability at each $k$.}
\end{figure}

The Fig.~\ref{fig:6} is showing cumulative degree distribution of $G^{(6)}$ with seed graph $G^{(0)}=K_3$. The exponential decay of degree distribution is also observed in \cite{jung2002geometric} when the derived resultant graph has same number of offsprings for all nodes. The proposed model in \cite{zhang2006deterministic} has also exponential degree distribution.

\subsection{Diameter}
Diameter of a graph is the longest shortest path in the graph. We show in the following theorem that as a Corona graph grows, the diameter also increases.
\begin{mydef-theorem}\label{thm:diam}
Consider the Corona graph $G^{(m)}$ generated by a seed graph $G^{(0)}$ with the diameter $D^{(0)}$. Then, the diameter $D^{(m)}$ of $G^{(m)}$ is $D^{(0)}+2m.$
\end{mydef-theorem}
\begin{proof}
We prove it by induction method. Let $i,j$ be the end nodes of the diameter $D^{(0)}$ in $G^{(0)}$. In $G^{1}=G^{(0)}\circ G^{(0)}$, a single instance of $G^{(0)}$ will be attached to every node(according to definition of corona product) including $i$ and $j$ and also each node will be connected to each node of $G^{(0)}$. Now, the diameter($D^{(1)}$) of the graph  will be elongated as $D^{(0)}+\mbox{(an edge from i)} + \mbox{(an edge from j)} = D^{(0)}+2$.

For inductive hypothesis, let the diameter of the Corona graph $G^{(m-1)}$ be $D^{(m-1)}=D^{(0)}+2(m-1)$ and $(i^{'},j^{'})$ are the two end nodes of $D^{(m-1)}$. Again, in $G^{(m)}$, a single instance of $G^{(0)}$ will be attached to every node of $G^{(m-1)}$ including $(i^{'},j^{'})$ and hence $D^{(m)}$ in $G^{(m)}$ is $D^{(0)}+2(m-1)+\mbox{(an edge from $i^{'}$)} + \mbox{(an edge from $j^{'}$)}$ $=D^{(0)}+2m$.
\end{proof}

\subsection{Betweenness Distribution}
Betweenness evaluates the importance of a person in a social network (see \cite{boccaletti2006complex} and the references therein). It is also being assessed in the models studying the cascading failures as the load on the $i^{th}$ node \cite{newman2003structure},\cite{holme2002edge},\cite{holme2002vertex}. Mathematically, it is formulated as
\begin{equation}\label{eqn:12}
b_i=\sum\limits_{j\ne i\ne k}\frac{\sigma_{jk}(i)}{\sigma_{jk}}
\end{equation}
where $b_i$ is the betweenness centrality of the $i^{th}$ node, $\sigma_{jk}(i)$ is the number of the shortest paths between nodes $j$ and $k$ through $i$, and $\sigma_{jk}$ is the total shortest paths between $j$ and $k$. 

It is difficult to find out all the pairs of shortest paths for $\sigma_{jk}(i)$ and $\sigma_{jk}$ for the Corona graphs generated by a seed graph. We observe that there is a distinct shortest path in a Corona graph between all pair of nodes when the seed graph is a clique. There are two reasons of this observation--(1) each node in a clique is connected to each other node of the clique and hence making a unique shortest path between each other. (2) each seed graph $G^{(0)}$ is connected in the $i^{th}$ corona product to any of the node $n_j$ of $G^{(i-1)}$ and this node $n_j$ is connected with a unique path to all nodes of this $G^{(0)}$ as well as it is acting as a unique bridge (and hence a unique path) between the nodes of $G^{(0)}$ and the rest of the nodes of $G^{(i-1)}$ and $G^{(i)}$. Hence, when $G^{(0)}=K_n$ (where $n\ge 2$), eqn.(13) can be modified as
\begin{equation}\label{eqn:13}
b_i=\sum\limits_{j\ne i\ne k}\sigma_{jk}(i)
\end{equation}
Now, for assessing the number of shortest paths, we use reasoning similar to \cite{qi2009structural},\cite{ghim2004packet}. Let all the nodes attached to node $i$ be termed as offsprings. Similarly, the nodes attached to the offsprings as well as the nodes attached recursively to all the offsprings in the subsequent corona products $G^{(\tau)}$ (corresponding to node $i$) are being termed as offsprings. 

\begin{mydef-theorem}\label{thm:betw}
The betweenness distribution of a Corona graph generated by a clique as the seed graph follows power law with exponent approximately $2.$
\end{mydef-theorem}
\begin{proof}
There are two types of shortest paths in the Corona graph $G^{(t)}$ generated by the seed graph $K_n$ through node $i$: 
\begin{enumerate}
\item[1.] A $=$ Total shortest paths between offsprings to the rest of the nodes of the graph.
\item[2.] B $=$ Total shortest paths between offsprings.
\end{enumerate}
such that $b_t^{i}(\tau)=A+B$ where $b_t^{i}(\tau)$ represents the total number of shortest paths through the node $i$ which was added in $\tau^{th}$ ($\tau\le t$) step during the formation of $G^{(t)}$. For $t\gg 1$,
\begin{equation*}
b_t(\tau)\approx A = |V^{(\tau)}|(|V^{(t)}|-|V^{(\tau)}|)
\end{equation*}
where $|V^{(\tau)}|$ and $|V^{(t)}|$ represent the total offsprings in $G^{(\tau)}$ (where $\tau\ge 1$) (corresponding to node $i$) and number of nodes in $G^{(t)}$ such that
\begin{eqnarray}
|V^{(\tau)}| &=& \sum\limits_{j=1}^{\tau}n(n+1)^{i-1}\nonumber\\
			 &=& ((n+1)^{\tau}-1).\label{eqn:14}
\end{eqnarray}
and
\begin{equation}\label{eqn:15}
|V^{(t)}|=n(n+1)^{t}.
\end{equation}
Using equations (~\ref{eqn:14}) and (~\ref{eqn:15}), we get
\begin{eqnarray}
b_t(\tau) &=& ((n+1)^{\tau}-1) (n(n+1)^{t}-n+1)^{\tau}+1).\nonumber\\
		\mbox{For $t\ggg 1$ and $\tau\gg 1$, we get}\nonumber\\
b_t(\tau) &\approx & n(n+1)^{t+\tau-1}.
\end{eqnarray}\label{eqn:16}
Let cumulative betweenness distribution be denoted by $P_c(b)$. Then,
\begin{figure}
\centering
     \includegraphics[width=0.5\textwidth]{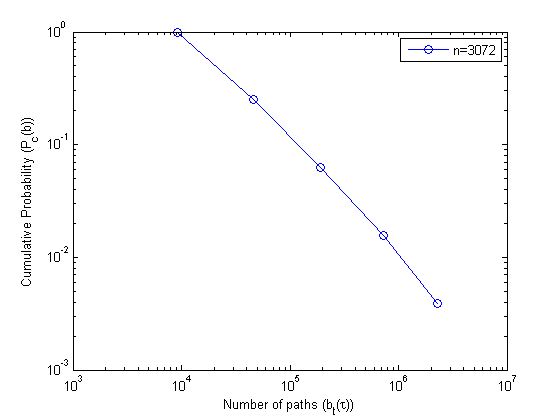} 
\captionsetup{justification=raggedright,singlelinecheck=false}
\caption{\label{fig:7}Cumulative power-law betweenness distribution of $G^{(5)}$ generated with seed graph $G^{(0)}=K_3$ as an instance of a clique. The betweenness values are generated with igraph \cite{csardi2006igraph}.}
\end{figure}
\begin{eqnarray}
P_c(b)&=&\sum\limits_{\mu=1}^{t-\tau+1}\frac{|V^{(\mu)}|}{|V^{(t)}|}\nonumber\\
      &\approx & \frac{n (n+1)^{t-\tau+1}}{n(n+1)^{t}}\nonumber\\
      &=& (n+1)^{1-\tau}\nonumber\\
      &\approx & \frac{|V^{(t)}|}{n(n+1)^{t+\tau-1}}\nonumber\\
      &\approx & (n(n+1)^{t+\tau-1})^{-1}.\nonumber
\end{eqnarray}
Thus
\begin{eqnarray}
      (n(n+1)^{t+\tau-1})^{1-\gamma_b}&=& (n(n+1)^{t+\tau-1})^{-1}\nonumber
\end{eqnarray}
which implies
\begin{equation}
      \gamma_b \approx 2.\label{eqn:17}
\end{equation}
\end{proof}
The Fig.~\ref{fig:7} is showing cumulative power-law betweenness distribution for $G^{(5)}$ generated by the clique $G^{(0)}=K_3$ as the seed graph. The betweenness distribution featuring power-law was also observed in \cite{goh2001universal},\cite{goh2003goh}, \cite{vazquez2002large}. In \cite{qi2009structural}, the authors analysed a family of deterministic recursive trees having exponential decay in cumulative degree distribution as well as cumulative power-law betweenness distribution.

\section{Spectra of Corona Graphs}\label{sec:Spec}
Let $G^{(m)}$ be the Corona graph generated by the seed graph $G$. The adjacency matrix $A(G^{(m)})$ associated with $G^{(m)}$ is given by
$$\mathbf{A(G^{(m)})} = \begin{bmatrix}
A(G^{(m-1)}) &  \mathbf{1}^{T}_{n(n+1)^{m-1}} \otimes I_{n(n+1)^{m-1}}\\
\mathbf{1}_{n(n+1)^{m-1}} \otimes I_{n(n+1)^{m-1}} & A(G) \otimes I_{n(n+1)^{m-1}} \\
\end{bmatrix}$$
\noindent where $A(G^{(m-1)})$ is the adjacency matrix of $G^{(m-1)}$, $I_{n(n+1)^{m-1}}$ is the identity matrix of order $n(n+1)^{m-1}$ and $\mathbf{1}_{n(n+1)^{m-1}}$ is the all-one vector of order $n(n+1)^{m-1}$.
\noindent We denote spectra of the corona graph $G^{(m)}$ i.e. spectrum of $A(G^{(m)})$ by
\begin{equation}
\sigma(G^{(m)})=\{\lambda_{1},\lambda_{2},\ldots,\lambda_{n(n+1)^m}\}
\end{equation} 
where $\lambda_{1}\le\lambda_{2}\le\ldots\le\lambda_{n(n+1)^m}$ and spectral radius of $A(G^{(m)})$ is denoted as $\rho(G^{(m)})$.
The spectrum of corona product of two graphs $G=G_{1}\circ G_{2}$ where $G_1$ is any graph and $G_{2}$ is a regular graph, and the Laplacian spectrum of corona product of any two graphs are provided by Barik et. al. in \cite{barik2007spectrum}. Inspired by their work, we derive $\sigma(G^{(m)})$ when the seed graph $G$ is regular. 

In the next theorem, we provide $\sigma(G^{(m)})$ in terms of the eigenvalues of the seed graph $G^{(0)}= G$. 

\begin{figure}
\centering

  \subfloat[]
  {%
    \adapttobigpicture
   \includegraphics[width=0.5\textwidth]{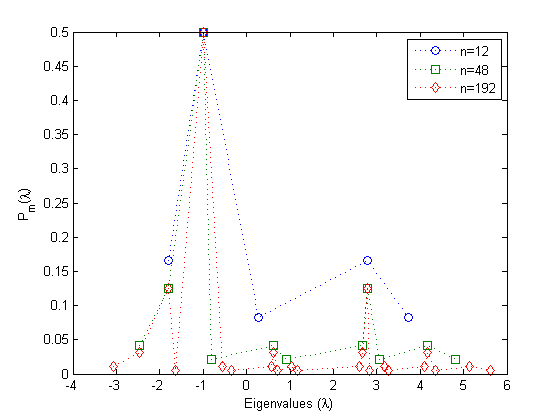} 
   \label{8a}
  }
  \subfloat[]
  {%
    \adapttobigpicture
   \includegraphics[width=0.5\textwidth]{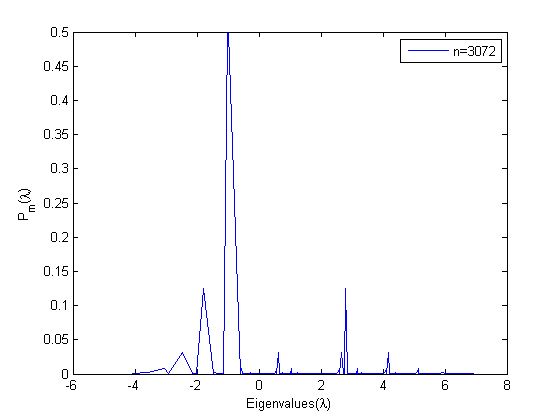} 
      \label{8b}
  }
\caption{ \label{8}Eigenvalues distribution with the seed graph $G^{(0)}=K_3$ for (a) $G^{(1)}$,$G^{(2)}$ and $G^{(3)}$ having 12,48 and 192 nodes respectively. (b) $G^{(5)}$ having 3072 nodes. Here, $P_m(\lambda)$ is showing the probability distribution for $G^{(m)}$.}
\end{figure}

\begin{mydef-theorem}\label{thm:spectra}
Let $G^{(0)}= G$ be a regular graph such that $\sigma(G)=(\mu_{1},\mu_{2},\ldots,\mu_{n}=r)$ (where, $\mu_{1}\le\mu_{2}\le\ldots\le\mu_{n}=r$) and spectral radius of $G$ be $\rho(G)$. Then, $\sigma(G^{(m)})$ is given by
\begin{itemize}
\item[(a)] $\lambda_{i}=\dfrac{\mu_{i}+r\sum_{a=0}^{m-1}2^{a}\pm \biggl(\sum_{c=1}^{m-1}z_{c}+\sqrt{((r-\mu_i)\pm \sum_{c=1}^{m-1}z_{c})^2+2^{2m}.n}\biggr)}{2^m} \in \sigma(G^{(m)})$, with multiplicity $1$ for $i=1,\ldots,n(n+1)^{m-1}$ \\
where, \\$z_{k-1}=z_{k-2}+\sqrt{((r-\mu_{i})\pm z_{k-2})^2 + n.2^{2(k-1)}}$, $2\le k\le m$, $z_0=0$.
\item[(b)] $\mu_{i}\in\sigma(G^{(m)})$, with multiplicity $n(n+1)^{m-1}$ for $i=1,\ldots,n-1$.
\end{itemize}
The spectral radius of $G^{(m)}$ is given by$$\rho(G^{(m)})=\frac{\mu_n+r\sum_{a=0}^{m-1}2^{a}+ \bigl(\sum_{c=1}^{m-1}z_{c}+\sqrt{((r-\mu_{i})- \sum_{c=1}^{m-1}z_{c})^2+2^{2m}.n}\bigr)}{2^m}$$ \\where $z_{c}$ is defined above.
\end{mydef-theorem}
\begin{proof}
Let $\mathbf{1}_{n}$ be the all-one vector of dimension $n$. Let $Z_1,Z_2,\ldots,Z_n$ be the eigenvectors associated with eigenvalues $\lambda_1,\lambda_2,\ldots,\lambda_n$ of $A(G)$ respectively. Then, the spectrum of $G^{(1)}=G^{(0)}\circ G^{(0)}$ (as given in Theorem 3.1 of \cite{barik2007spectrum}) is as follows
\begin{itemize}
\item[(i)] $\frac{\mu_i+r\pm\sqrt{(r-\mu_i)^2+4n}}{2}\in$ $\sigma(G^{(1)})$ with multiplicity $1$ for $i=1,\dots,n.$ Let $\lambda_i$ and $\hat{\lambda}_i$ be their representations. Then, eigenvectors corresponding to them are
$$ 
\begin{pmatrix}
Z_i\\
\frac{1}{\lambda_i-r}(\mathbf{1}_{n}\otimes Z_i)
\end{pmatrix}
,
\begin{pmatrix}
Z_i\\
\frac{1}{\hat{\lambda}_i-r}(\mathbf{1}_{n}\otimes Z_i)
\end{pmatrix}
$$
\item[(ii)] $\mu_i\in\sigma(G^{(1)})$ with multiplicity $n$ for $i=1,\dots,n-1$ and the corresponding eigenvector is 
$$
\begin{pmatrix}
\mathbf{0}\\
Z_i\otimes e_j
\end{pmatrix}
$$
where $j\in[1,n]$ and $e_1,e_2,\ldots,e_n$ are the unit vectors of standard basis for $n\times n.$
\end{itemize}
Now, the eigenvalues for $G^{(m)}=G^{(m-1)}\circ G$ are
\begin{itemize}
\item[(a)] $\lambda=\dfrac{\mu_{i}+r\sum_{a=0}^{m-1}2^{a}\pm \biggl(\sum_{c=1}^{m-1}z_{c}+\sqrt{((r-\mu_i)\pm \sum_{c=1}^{m-1}z_{c})^2+2^{2m}.n}\biggr)}{2^m} \in \sigma(G^{(m)})$, with multiplicity $1$ for $i=1,\ldots,n(n+1)^{m-1}$ corresponding to the eigenvectors
$$
U^{(m)}_{l_m}=\begin{pmatrix}
			U^{(m-1)}_{l_{m-1}}\\
			\frac{1}{\lambda_{l_m}-r}\begin{pmatrix}
										\mathbf{1}_n\otimes U^{(m-1)}_{l_{m-1}}
									 \end{pmatrix}
		\end{pmatrix}
$$
where, 
$$
U^{(m-1)}_{l_{m-1}}=\begin{pmatrix}
			U^{(m-2)}_{l_{m-2}}\\
			\frac{1}{\lambda_{l_{m-1}}-r}\begin{pmatrix}
											\mathbf{1}_n\otimes U^{(m-2)}_{l_{m-2}}
										 \end{pmatrix}
		  \end{pmatrix},
		  \ldots\,,
U^{(1)}_{l_1}=\begin{pmatrix}
			Z_i\\
			\frac{1}{\lambda_{l_{1}}-r}\begin{pmatrix}
											\mathbf{1}_n\otimes Z_i
										 \end{pmatrix}
		\end{pmatrix}.
$$ 
$l_a\in[1,2n(n+1)^{a-1}]$ such that $\lambda_{l_a}$ is an eigenvalue of $G^{(a)}$, 
$z_{k-1}=z_{k-2}+\sqrt{((r-\mu_{i})\pm z_{k-2})^2 + n.2^{2(k-1)}}$, $2\le k\le m$, $z_0=0$.
\item[(b)] $\mu_{i}\in\sigma(G^{(m)})$, with multiplicity $n(n+1)^{m-1}$ for $i=1,\ldots,n-1$ and the eigenvector corresponding to these eigenvalues is 
$$
\begin{pmatrix}
\mathbf{0}\\
Z_i\otimes e_j
\end{pmatrix}
$$
where $j\in[1,n(n+1)^{m-1}]$ and $e_j$ are the unit vectors of standard basis for $n(n+1)^{m-1}\times n(n+1)^{m-1}.$
\end{itemize}
\end{proof}
Fig.~\ref{8} is showing the eigenvalue distribution of Corona graphs, with the seed graph $G^{(0)}=K_3$. It is the ratio between the multiplicity of eigenvalues and the number of nodes. 
We derive the spectra of star graphs $S_k$ with $k\ge 3$ which is a special case of irregular graph.
\begin{mydef-theorem}\label{thm:tree_spectra}
Let $k\ge$ $3$ be an integer. The spectrum of the graph $S_k \circ S_k$ consists of the following eigenvalues
\begin{itemize}
\item[(a)] $\lambda_{z}=X_1^z + \frac{\mu_i}{3}\in\sigma(G^{(1)})$ with multiplicity $1$, and
\item[(b)] $0\in\sigma(G^{(1)})$ with multiplicity $k(k-2)$. 
\end{itemize}
where, $\lambda_{z}$ are the eigenvalues for $G^{(1)}$ with $z=1,2,3$ for the $3$ angles i.e. $\frac{\theta}{3},\frac{2\pi+\theta}{3},\frac{4\pi+\theta}{3}$ as shown in following sub-expressions\\ 
$X^z=\frac{2}{3}\cos\frac{y\pi+\theta}{3}\sqrt{\mu_i^2+(6k-3)}$, $\theta=\cos^{-1}\biggl(\frac{2\mu_i^3+\mu_i(18-9k)+(54k-54)}{2(\mu_i^2+(6k-3))^\frac{3}{2}}\biggr)-y\pi$\\
where $y=0,2,4$. Here, $\lambda_{z}\in\Bigl[-\mu_i+\frac{2\mu_i}{3}-\frac{2\sqrt{(6k-3)}}{3},\mu_i+\frac{2\sqrt{(6k-3)}}{3}\Bigr].$
\end{mydef-theorem}
\begin{proof}
Let $\mathbf{1}_{n}$ be the all-one vector of dimension $n$. Let $Z_1,Z_2,\dots,Z_n$ be the eigenvectors corresponding to seed graph's eigenvalues $\mu_1,\mu_2,\dots,\mu_n$ respectively. The eigenvalues corresponding to $G^{(1)}$ are $\lambda_{z}=X_1^z + \frac{\mu_i}{3}\in\sigma(G^{(1)})$ for $z=1,2,3$ and the eigenvectors corresponding to them are as
\begin{equation*}
\begin{pmatrix}
Z_i\\
\frac{\lambda_i+1}{\lambda_i^2-k+1}(\mathbf{1}_{n}\otimes Z_i)
\end{pmatrix}
\end{equation*}
$G^{(1)}$ has $0$ as its one of the eigenvalues and its  multiplicity is of $k(k-2)$.
\end{proof}

\begin{figure}
\centering

  \subfloat[]
  {%
    \adapttobigpicture
   \includegraphics[width=0.5\textwidth]{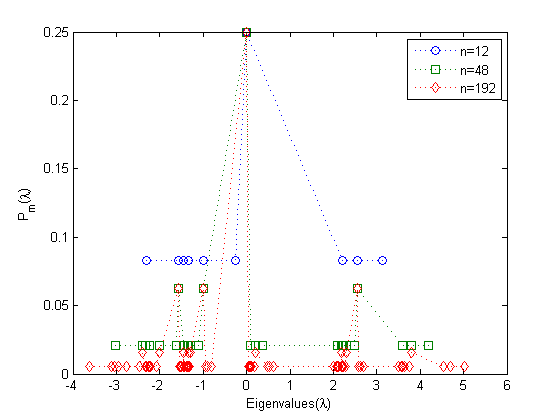} 
   \label{9a}
  }
  \subfloat[]
  {%
    \adapttobigpicture
   \includegraphics[width=0.5\textwidth]{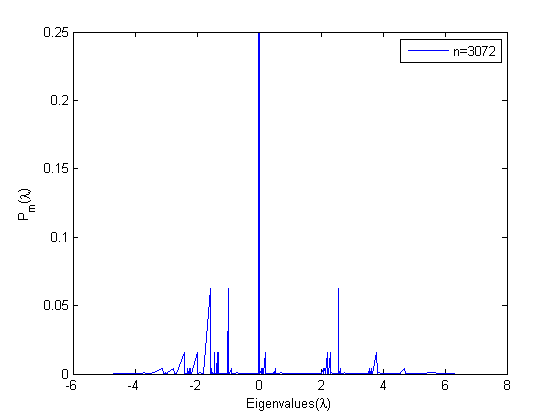} 
      \label{9b}
  }
\caption{ \label{9}Eigenvalues distribution with the seed graph $G^{(0)}=S_3$ for (a) $G^{(1)}$,$G^{(2)}$ and $G^{(3)}$ having 12,48 and 192 nodes respectively. (b) $G^{(5)}$ having 3072 nodes. Here, $P_m(\lambda)$ is showing the probability distribution for $G^{(m)}$.}
\end{figure}

\begin{mydef-corollary}
Let $G$ be the seed star graph $S_k$ for each $k\ge 3$ such that $\sigma(G)=\{\mu_1,\mu_2,\dots,\mu_n\}$. Let $m\ge 1$.Then $\sigma(G^{(m)})$ is given by
\begin{itemize}
\item[(a)] $\lambda_{z,1}=X_1^z + \frac{\mu_i}{3}$,
\ldots,
$\lambda_{z,m}=\sum_{j=0}^{m-1}(\frac{1}{3})^{j} X_{m-j}^z + (\frac{1}{3})^{m-1}(\frac{\mu_i}{3})\in\sigma(G^{(m)})$ with multiplicity $1$ for each of them,
\item[(b)] $0\in\sigma(G^{(m)})$ with multiplicity $k(k-2)(k+1)^{(m-1)}$.
\end{itemize}
where, $\lambda_{z,j}$ are the eigenvalues for $G^{(m)}$ such that $j$ represents $j^{th}$ corona product with $z=1,2,3$ for the $3$ angles i.e. $\frac{\theta_j}{3},\frac{2\pi+\theta_j}{3},\frac{4\pi+\theta_j}{3}$ as shown in following sub-expressions\\ $X_l^z=\frac{2}{3}\cos\frac{y\pi+\theta_l}{3}\sqrt{\bigl(\sum_{j=1}^{l-1}(\frac{1}{3})^{j-1}X_{l-j}^z+(\frac{1}{3})^{l-1}\mu_i\bigr)^2+(6k-3)}$,\\
$X_1^z=\frac{2}{3}\cos\frac{y\pi+\theta_1}{3}\sqrt{\mu_i^2+(6k-3)}$\\
$\theta_{m}=\cos^{-1}\biggl(\frac{2(\sum_{j=0}^{m-1}(\frac{1}{3})^{j} X_{m-j}^z + (\frac{1}{3})^{m-1}(\frac{\mu_i}{3}))^3-(9k-18)(\sum_{j=0}^{m-1}(\frac{1}{3})^{j} X_{m-j}^z + (\frac{1}{3})^{m-1}(\frac{\mu_i}{3}))+(54k-54)}{2((\sum_{j=0}^{m-1}(\frac{1}{3})^{j} X_{m-j}^z + (\frac{1}{3})^{m-1}(\frac{\mu_i}{3}))^2+(6k-3))^\frac{3}{2}}\biggr)-y\pi$,\\
where $m$ is the $m^{th}$ corona product and $y=0,2,4$ for the three angles.\\
Here, $\lambda_{z,j}\in\Bigl[-\mu_i+\frac{2\mu_i}{3^j}-\frac{2j\sqrt{(6k-3)}}{3},\mu_i+\frac{2j\sqrt{(6k-3)}}{3}\Bigr]$,
where $j=1,\dots,m$.
\end{mydef-corollary}
\begin{proof}
Let $\mathbf{1}_{n}$ be the all $1$ vector of dimension $n$. The proof follows by using similar arguments given in the proof of Theorem ~\ref{thm:spectra} and Theorem ~\ref{thm:tree_spectra}. The eigenvectors corresponding to the eigenvalues $\lambda_{z,j}$, $z=1,2,3$,  $j\in[1,m]$ are given by
$$
U^{(l)}=\begin{pmatrix}
			U^{(l-1)}\\
			\frac{\lambda_{z,j_l}+k-1}{\lambda_{z,j_l}^2-k+1}\begin{pmatrix}
																\mathbf{1}_n\otimes U^{(l-1)}
															 \end{pmatrix}
		\end{pmatrix}
$$
where 
$$
U^{(l-1)}=\begin{pmatrix}
			U^{(l-2)}\\
			\frac{\lambda_{z,j_{l-1}}+k-1}{\lambda_{z,j_{l-1}}^2-k+1}\begin{pmatrix}
																		\mathbf{1}_n\otimes U^{(l-2)}
										 							 \end{pmatrix}
		  \end{pmatrix},
		  \ldots\,,
U^{(1)}=\begin{pmatrix}
			Z_i\\
			\frac{\lambda_{z,j_1}+k-1}{\lambda_{z,j_1}^2-k+1}\begin{pmatrix}
																\mathbf{1}_n\otimes Z_i
															 \end{pmatrix}
		\end{pmatrix},
$$
$\lambda_{z,j_l}$ are the eigenvalues for $G^{(l)}$, $1\le l\le m$.

$G^{(m)}$ has 0 as its one of the eigenvalue and its  multiplicity is of $k(k-2)(k+1)^{(m-1)}$.

\end{proof}
In Fig.~\ref{9}, we plot the distribution of eigenvalues generated with the seed graph $G^{(0)}=S_3$.

\noindent The Laplacian matrix $L(G^{(m)})$ associated with $G^{(m)}$ has the form
$$\mathbf{L(G^{(m)})} = \begin{bmatrix}
\mathbf{L(G^{(m-1)})}+n\mathbf{I}_{n(n+1)^{m-1}} &  -\mathbf{1}^{T}_{n(n+1)^{m-1}} \otimes \mathbf{I_{n(n+1)^{m-1}}}\\
-\mathbf{1}_{n(n+1)^{m-1}} \otimes \mathbf{I_{n(n+1)^{m-1}}} & \mathbf{(L(G)}+\mathbf{I}_{n}) \otimes \mathbf{I_{n(n+1)^{m-1}}} \\
\end{bmatrix}$$

\noindent where $L(G^{(m-1)})$ is the Laplacian matrix of $G^{(m-1)}$, $I_{n}$ and $I_{n(n+1)^{m-1}}$ are the identity matrices.
\noindent We denote the Laplacian spectra $L(G^{(m)})$ of $G^{(m)}$ by
\begin{equation}
S(G^{(m)})=\{\lambda_{1},\lambda_{2},\ldots,\lambda_{n(n+1)^{(m)}}\}
\end{equation} 
where $0=\lambda_1\le\lambda_2\le\ldots\le\lambda_{n(n+1)^{(m)}}$. In the following theorem, we determine the elements of $S(G^{(m)})$ in terms of the Laplacian eigenvalues of the seed graph, $G^{(0)}= G$ where $S(G)=\{0=\nu_{1},\nu_{2},\ldots,\nu_{n}\}$ (where, $\nu_{1}\le\nu_{2}\le\ldots\le\nu_{n}$) is the Laplacian spectra of $G$. The algebraic connectivity of a graph is defined as the second smallest eigenvalue of $L(G)$ \cite{bapat2010graphs}.

\begin{figure}
\centering

  \subfloat[]
  {%
    \adapttobigpicture
   \includegraphics[width=0.5\textwidth]{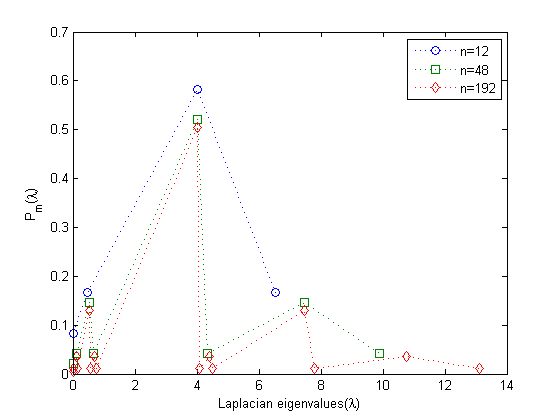} 
   \label{10a}
  }
  \subfloat[]
  {%
    \adapttobigpicture
   \includegraphics[width=0.5\textwidth]{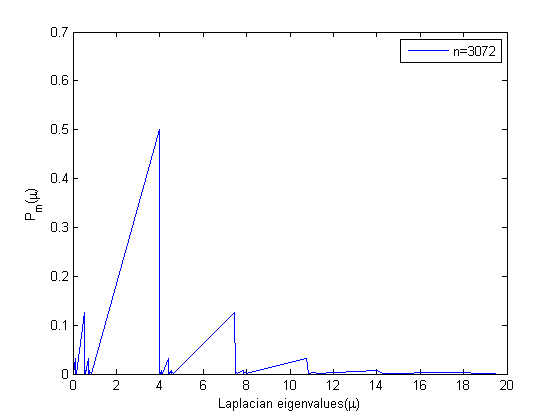} 
      \label{10b}
  }
\caption{ \label{10}Laplacian eigenvalues distribution with the seed graph $G^{(0)}=K_3$ for (a) $G^{(1)}$,$G^{(2)}$ and $G^{(3)}$ having 12,48 and 192 nodes respectively. (b) $G^{(5)}$ having 3072 nodes. Here, $P_m(\lambda)$ is showing the probability distribution for $G^{(m)}$.}
\end{figure}

\begin{mydef-theorem}\label{thm:laplacian}
Let $G$ be a simple connected graph. We denote the algebraic connectivity of $G$ and $G^{(m)}$ by $a(\nu_2)$ and $a(\lambda_2)$ respectively. Then, Laplacian spectra $S(G^{(m)})$ of $G$ is given by
\begin{itemize}
\item[(a)] $\dfrac{\nu_{i}+(n+1)\sum_{i=0}^{m-1}2^{i}\pm\sum_{i=1}^{m}z_{i}}{2^m}$ $\in S(G^{(m)})$  with multiplicity $1$ for $i=1,\ldots,n(n+1)^{m-1}$.
where, \\$z_{1}=\sqrt{(\nu_i+n+1)^2-4\nu_{i}}$ ,\\
 $$z_{k}=\sqrt{(\nu_i+(n+1)\sum_{i=0}^{k-1}2^i\pm\sum_{i=1}^{k-1}z_i)^2-2^{(k+1)}(\nu_i+(n+1)\sum_{i=0}^{k-2}2^i\pm\sum_{i=1}^{k-1}z_i)}$$
 for $2\le k\le m.$
\item[(b)] $\nu_{i}+1\in S(G^{(m)})$ with multiplicity $n(n+1)^{m-1}$ for $i=2,\ldots,n$.
\end{itemize}
Hence, the algebraic connectivity of $S(G^{(m)})$ is $$a(\lambda_2)=\dfrac{\nu_{2}+(n+1)\sum_{i=0}^{m-1}2^{i}-\sum_{i=1}^{m}z_{i}}{2^m}< 1 $$ where $z_i$ can be defined as above.
\end{mydef-theorem}
\begin{proof}
Let $\mathbf{1}_{n}$ be the all-one vector of dimension $n$. Let $X_1,X_2,\ldots,X_n$ be the eigenvectors of corresponding to the eigenvalues $\nu_1,\nu_2,\ldots,nu_n$ respectively of $L(G)$. Then, the spectrum of $G^{(1)}=G^{(0)}\circ G^{(0)}$ (as defined in Theorem 3.2 of \cite{barik2007spectrum}) is as follows
\begin{itemize}
\item[(i)] $\frac{\nu_i+n+1\pm\sqrt{(\nu_i+n+1)^2-4\nu_i}}{2}\in$ $S(G^{(1)})$ with multiplicity $1$ for $i=1,\dots,n$. Let $\mu_i$ and $\hat{\mu}_i$ be their representation. Then, eigenvector corresponding to them is as
$$ 
\begin{pmatrix}
X_i\\
\frac{1}{1-\mu_i}(\mathbf{1}_{n}\otimes X_i)
\end{pmatrix}
,
\begin{pmatrix}
X_i\\
\frac{1}{1-\hat{\mu}_i}(\mathbf{1}_{n}\otimes X_i)
\end{pmatrix}
$$
\item[(ii)] $\nu_i+1\in S(G^{(1)})$ with multiplicity $n$ for $i=2,\dots,n$ and the corresponding eigenvector is 
$$
\begin{pmatrix}
\mathbf{0}\\
X_i\otimes e_j
\end{pmatrix}
$$
where $j\in[1,n]$ and $e_1,e_2,\ldots,e_n$ are the unit vectors of standard basis for $n\times n.$
\end{itemize}
Now, the laplacian eigenvalues for $G^{(m)}=G^{(m-1)}\circ G$ is
\begin{itemize}
\item[(a)] $\dfrac{\nu_{i}+(n+1)\sum_{i=0}^{m-1}2^{i}\pm\sum_{i=1}^{m}z_{i}}{2^m}$ $\in S(G^{(m)})$  with multiplicity $1$ for $i=1,\ldots,n(n+1)^{m-1}$.
\\ corresponding to the eigenvectors
 $$
 U^{(m)}_{l_m}=\begin{pmatrix}
 			U^{(m-1)}_{l_{m-1}}\\
 			\frac{1}{\lambda_{l_m}-r}\begin{pmatrix}
 										\mathbf{1}_n\otimes U^{(m-1)}_{l_{m-1}}
 									 \end{pmatrix}
 		\end{pmatrix}
 $$
 where, 
 $$
 U^{(m-1)}_{l_{m-1}}=\begin{pmatrix}
 			U^{(m-2)}_{l_{m-2}}\\
 			\frac{1}{\lambda_{l_{m-1}}-r}\begin{pmatrix}
 											\mathbf{1}_n\otimes U^{(m-2)}_{l_{m-2}}
 										 \end{pmatrix}
 		  \end{pmatrix},
 		  \ldots\,,
 U^{(1)}_{l_1}=\begin{pmatrix}
 			Z_i\\
 			\frac{1}{\lambda_{l_{1}}-r}\begin{pmatrix}
 											\mathbf{1}_n\otimes Z_i
 										 \end{pmatrix}
 		\end{pmatrix}.
 $$ 
 $l_a\in[1,2n(n+1)^{a-1}]$ such that $\Delta_{l_a}$ is an eigenvalue of $G^{(a)}$, \\$z_{1}=\sqrt{(\nu_i+n+1)^2-4\nu_{i}}$ ,
  $$z_{m}=\sqrt{(\nu_i+(n+1)\sum_{i=0}^{m-1}2^i\pm\sum_{i=1}^{m-1}z_i)^2-2^{(m+1)}(\nu_i+(n+1)\sum_{i=0}^{m-2}2^i\pm\sum_{i=1}^{m-1}z_i)}$$
  for $m\ge 2.$
\item[(b)] $\nu_{i}+1\in S(G^{(m)})$, with multiplicity $n(n+1)^{m-1}$ for $i=2,\ldots,n$ and the eigenvector corresponding to these eigenvalues is 
$$
\begin{pmatrix}
\mathbf{0}\\
X_i\otimes e_j
\end{pmatrix}
$$
where $j\in[1,n(n+1)^{m-1}]$ and $e_j$ are the unit vectors of standard basis for $n(n+1)^{m-1}\times n(n+1)^{m-1}.$
\end{itemize}

\end{proof}
Fig.~\ref{10} is showing the laplacian eigenvalue distribution, with the seed graph $G^{(0)}=K_3$. 

\noindent The signless Laplacian matrix $S_Q(G^{(m)})$ of  $G^{(m)}$ is of the form
$$\mathbf{S_Q(G^{(m)})} = \begin{bmatrix}
\mathbf{S_Q(G^{(m-1)})}+n\mathbf{I}_{n(n+1)^{m-1}} &  \mathbf{1}^{T}_{n(n+1)^{m-1}} \otimes \mathbf{I_{n(n+1)^{m-1}}}\\
\mathbf{1}_{n(n+1)^{m-1}} \otimes \mathbf{I_{n(n+1)^{m-1}}} & \mathbf{(S_Q(G)}+\mathbf{I}_{n}) \otimes \mathbf{I_{n(n+1)^{m-1}}} \\
\end{bmatrix}$$
\noindent where $S_Q(G^{(m-1)})$ is the signless Laplacian matrix of $G^{(m-1)}$, $I_{n}$ and $I_{n(n+1)^{m-1}}$ are the identity matrices.
\noindent The spectrum of signless Laplacian of $Q(G^{(m)})$ for $G^{(m)}$ is denoted by
\begin{equation}
S_Q(G^{(m)})=\{\lambda_1,\lambda_2,\ldots,\lambda_{n(n+1)^{(m)}}\}
\end{equation}
where $\lambda_1\le \lambda_2\le\ldots\le \lambda_{n(n+1)^{m}}$. Recently, there is a lot of work on signless Laplacian matrices and their $Q-$spectra as in \cite{cui2012spectrum},\cite{cvetkovic2009towards},\cite{haemers2004enumeration} and their authors think that it is more useful and hence a lot of work is going on to find their usefulness as \cite{cvetkovic2010towards},\cite{cvetkovic2010toward}. We will define the $S_Q(G^{(m)})$ inspired by the work of \cite{cui2012spectrum} on two graphs and here taking the seed graph $G$ as the $r$-regular graph. In the following theorem, we derive the elements of $S_Q(G^{(i)})$ in terms of the signless Laplacian eigenvalues of a regular seed graph, $G^{(0)}=G$ such that $S_Q(G)=(q_{1},q_{2},\ldots,q_{n}=2r)$ (where, $q_{1}\le q_{2}\le\ldots\le q_{n}=2r$). 



\begin{figure}
\centering

  \subfloat[]
  {%
    \adapttobigpicture
   \includegraphics[width=0.5\textwidth]{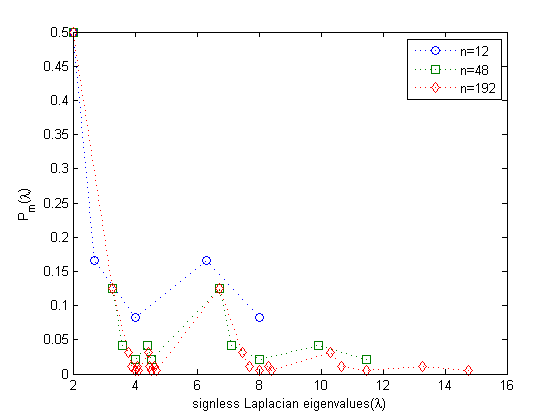} 
   \label{11a}
  }
  \subfloat[]
  {%
    \adapttobigpicture
   \includegraphics[width=0.5\textwidth]{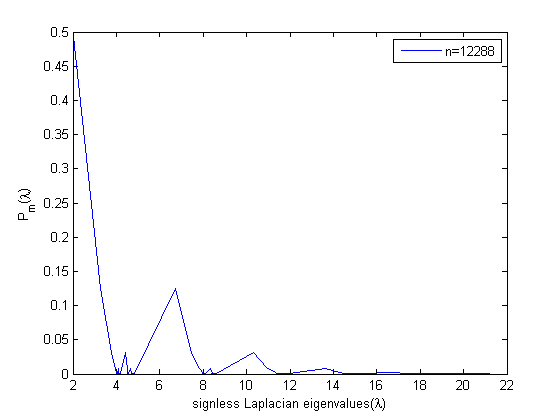} 
      \label{11b}
  }
\caption{ \label{11}Signless laplacian eigenvalues distribution with the seed graph $G^{(0)}=K_3$ for (a) $G^{(1)}$,$G^{(2)}$ and $G^{(3)}$ having 12,48 and 192 nodes respectively. (b) $G^{(5)}$ having 3072 nodes. Here, $P_m(\lambda)$ is showing the probability distribution for $G^{(m)}$.}
\end{figure}

\begin{mydef-theorem}\label{thm:signless_regular}
Let $G$ be a simple connected graph. Then, $S_Q(G^{(m)})$ is given by
\begin{itemize}
\item[(a)] $\lambda_i=\dfrac{q_i+n\sum_{i=0}^{m-1}2^i+r\sum_{i=1}^{m}2^i+\sum_{i=0}^{m-1}2^i\pm\sum_{j=1}^{m}z_j}{2^m}$ with multiplicity of $1$, for $i=1,\ldots,n(n+1)^{m-1}$ \\
where,\\ $$z_{j}=\sqrt{(q_i+n\sum_{i=0}^{j-1}2^i+r(\sum_{i=1}^{j-1}2^i-2^j)+(\sum_{i=0}^{j-2}2^i-2^{j-1})\pm\sum_{i=1}^{j-1}z_i)^2+2^{2j}.n}$$ for $j=2,\ldots,m$ and $z_1=\sqrt{((q_i+n)-(2r+1))^2+4n}$. 
\item[(b)] $q_j+1$ with the multiplicity of $n(n+1)^{m-1}$ for $j=1,\ldots,n-1$
\end{itemize}
Hence, spectral radius of $S_Q(G^{(m)})$ is\\
$$q(S_Q(G^{(m)}))=\dfrac{q_i+n\sum_{i=0}^{m-1}2^i+r\sum_{i=1}^{m}2^i+\sum_{i=0}^{m-1}2^i+\sum_{j=1}^{m}z_j}{2^m}$$
where $z_j$ is defined as above.
\end{mydef-theorem}
\begin{proof}
We can prove the Part (a) of the theorem by induction and hence, the $q(S_Q(G^{(m)}))$ will be  followed from the proof which is as follows \\
Base case: For $j=1$, $G^{(1)}=G^{(0)}\circ G$, the $S_Q(G^{(1)})$ can be defined according to Theorem 3.2 of \cite{cui2012spectrum} as
\begin{itemize}
\item[(i)] $\frac{q_i+n+2r+1\pm\sqrt{((q_i+n)-(2r+1))^2+4n}}{2}\in$ $S_Q(G^{(1)})$ with multiplicity $1$ for $i=1,\dots,n$
\item[(ii)] $q_i+1\in S_Q(G^{(1)})$ with multiplicity $n$ for $i=1,\dots,n-1$ 
\end{itemize}
Inductive hypothesis: Let $j=m-1$, $G^{(m-1)}=G^{(m-2)}\circ G$, the $S_Q(G^{(m-1)})$ can be defined as
\begin{itemize}
\item[(a)] $\lambda_i=\dfrac{q_i+n\sum_{i=0}^{m-2}2^i+r\sum_{i=1}^{m-1}2^i+\sum_{i=0}^{m-2}2^i\pm\sum_{i=1}^{m-1}z_i}{2^{m-1}}$ with multiplicity of $1$, for $i=1,\ldots,n(n+1)^{m-2}$ \\
where,\\ $$z_{j}=\sqrt{(q_i+n\sum_{i=0}^{j-1}2^i+r(\sum_{i=1}^{j-1}2^i-2^j)+(\sum_{i=0}^{j-2}2^i-2^{j-1})\pm\sum_{i=1}^{j-1}z_i)^2+2^{2j}.n}$$ for $j=2,\ldots,m-1$ and $z_1=\sqrt{((q_i+n)-(2r+1))^2+4n}$. 
\item[(b)] $q_i+1$ with the multiplicity of $n(n+1)^{m-2}$ for $i=1,\ldots,n-1$
\end{itemize}
Inductive step: For $j=m$, $G^{(m)}=G^{(m-1)}\circ G$, the $S_Q(G^{(m)})$ can be obtained by substituting the eigenvalues of $S_Q(G^{(m-1)})$ in place of $q_i$ of the base step and we will get the eigenvalues as stated in theorem.

The other eigenvalues are $q_i+1$  with the multiplicity of $n(n+1)^{m-1}$ for $j=1,\ldots,n-1$.
\end{proof}
Fig.~\ref{11} is showing the signless laplacian eigenvalue distribution, with the seed graph $G^{(0)}=K_3$. 

Consider star graph $S_k$ which is an irregular graph. In the theorem below, we determine explicit formula of signless Laplacian elements of $G^{(m)}=S_k^{(m)}$.

\begin{mydef-theorem}\label{thm:signless_tree_spectra}
Let $S_Q(G)=\{q_1,q_2,\dots,q_n\}$ with $q_1\le q_2\le \ldots\le q_n$. Then, $S_Q(G^{(1)})$ is given by
\begin{itemize}
\item[(a)] $\lambda_{z}=X_1^z + \frac{q_i+2k+2}{3}\in S_{Q}(G^{(1)})$ with multiplicity $1$, and
\item[(b)] $q_j+1\in S_{Q}(G^{(1)})$ where for $j=2,\dots,n-1$ with multiplicity of each of them as $k$
\end{itemize}
where, $\lambda_{z}$ are the eigenvalues for $G^{(1)}$ with $z=1,2,3$ for the $3$ angles i.e. $\frac{\theta}{3},\frac{2\pi+\theta}{3},\frac{4\pi+\theta}{3}$ as shown in following sub-expressions\\ 
$X^z=\frac{2}{3}\cos\frac{y\pi+\theta}{3}\sqrt{q_i^2+q_i(k-2)+(k+1)^2}$,\\ $\theta=\cos^{-1}\biggl(\frac{2q_i^3+(3k-6)q_i^2-3(k^2-k-2)q_i+(70k-94-12\sum_{a=1}^{k-2}(a+2)(k-a-1)}{2(q_i^2+q_i(k-2)+(k+1)^2)^{\frac{3}{2}}}\biggr)-y\pi$\\
where $y=0,2,4$. Here, 
$\lambda_z\in\biggl[ -q_i+\frac{2q_i}{3}+\frac{Y-(A+\sqrt{4B-A^2})}{3},q_i+\frac{Y+A+\sqrt{4B-A^2}}{3}\biggr] $\\
where $A=(k-2),B=(k+1)^2,Y=(2k+2).$
\end{mydef-theorem}
\begin{proof}
Let $\mathbf{1}_{n}$ be the all $1$ vector of dimension $n\times 1$. Let $Y_1,\dots,Y_n$ be the eigenvectors corresponding to seed graph's eigenvalues $q_1,\dots,q_n$.  The eigenvalues for $G^{(1)}$ are $\lambda_{z}=X_1^z + \frac{q_i+2k+2}{3}\in S_{Q}(G^{(1)})$ for $z=1,2,3$ and the eigenvectors corresponding to them are as
\begin{equation*}
\begin{pmatrix}
Y_i\\
\frac{q_i-k+1}{q_i^2-q_i(k+2)+(k+1)}(\mathbf{1}_{n}\otimes Y_i)
\end{pmatrix}
\end{equation*}
$q_j+1\in S_{Q}(G^{(1)})$ are the other eigenvalues for $j=2,\dots,n-1$ with multiplicity of each of them as $k$.
\end{proof}

\begin{figure}
\centering

  \subfloat[]
  {%
    \adapttobigpicture
   \includegraphics[width=0.5\textwidth]{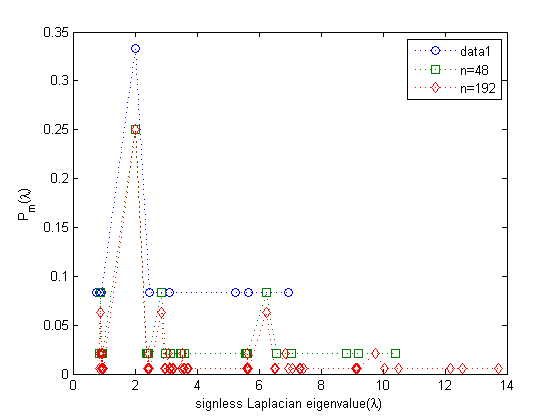} 
   \label{12a}
  }
  \subfloat[]
  {%
    \adapttobigpicture
   \includegraphics[width=0.5\textwidth]{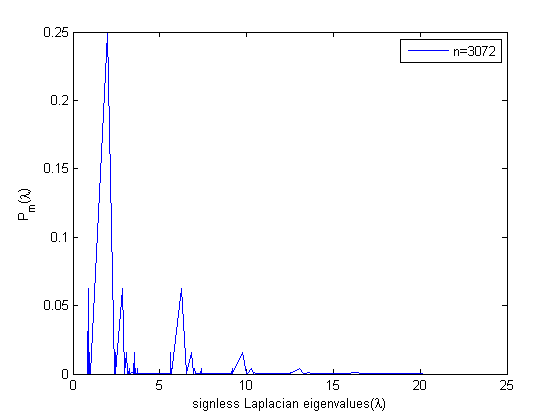} 
      \label{12b}
  }
\caption{ \label{12}Signless laplacian eigenvalues distribution with the seed graph $G^{(0)}=S_3$ for (a) $G^{(1)}$,$G^{(2)}$ and $G^{(3)}$ having 12,48 and 192 nodes respectively. (b) $G^{(5)}$ having 3072 nodes. Here, $P_m(\lambda)$ is showing the probability distribution for $G^{(m)}$.}
\end{figure}

\begin{mydef-corollary}
Let $S_Q(G)$ for star graph $S_k$ for each $k\ge 3$ is $\{q_1,q_2,\dots,q_n\}$ with $q_1\le q_2\le \ldots\le q_n$. Then, $S_Q(G^{(m)})$ is given by
\begin{itemize}
\item[(a)] $\lambda_{z,1}=X_1^z + \frac{q_i+2k+2}{3},\ldots,\lambda_{z,m}=\sum_{j=0}^{m-1}(\frac{1}{3})^j X_{m-j}+(\frac{1}{3})^m q_i+(2k+2)\sum_{j=1}^{m}(\frac{1}{3})^j\in S_{Q}(G^{(m)})$ with multiplicity $1$, and
\item[(b)] $q_i+1\in S_{Q}(G^{(m)})$ where for $i=2,\dots,n-1$ with multiplicity of each of them as $k(k+1)^{(m-1)}$.
\end{itemize}
where, $\lambda_{z,j}$ are the eigenvalues for $G^{(m)}$ (for $j^{th}$ corona product) with $z=1,2,3$ for the $3$ angles i.e. $\frac{\theta_j}{3},\frac{2\pi+\theta_j}{3},\frac{4\pi+\theta_j}{3}$ as shown in following sub-expressions\\ 
$X^z_l=\sqrt{(A)^2+A+(k+1)^2 }$\\
$X^z_1=\frac{2}{3}\cos\frac{y\pi+\theta}{3}\sqrt{q_i^2+q_i(k-2)+(k+1)^2}$,\\
where $A=(\sum_{j=0}^{m-2}(\frac{1}{3})^j X_{m-j-1}+(\frac{1}{3})^{m-1} q_i+(2k+2)\sum_{j=1}^{m-1}(\frac{1}{3})^j)$\\ $\theta=\cos^{-1}\biggl(\frac{2\lambda_{z,m-1}^3+(3k-6)\lambda_{z,m-1}^2-3(k^2-k-2)\lambda_{z,m-1}+(70k-94-12\sum_{a=1}^{k-2}(a+2)(k-a-1)}{2(\lambda_{z,m-1}^2+\lambda_{z,m-1}(k-2)+(k+1)^2)^{\frac{3}{2}}}\biggr)-y\pi$\\
where $y=0,2,4$ and $\lambda_{z,m-1}$ is as defined in part(a) of corollary. Here,\\
$\lambda_{z,1}\in\biggl[ -q_i+\frac{2q_i}{3}+\frac{Y-(A+\sqrt{4B-A^2})}{3},q_i+\frac{Y+A+\sqrt{4B-A^2}}{3}\biggr] $,
\ldots,
$\lambda_{z,m}\in\biggl[ -q_i+\frac{2q_i}{3^m}-\frac{m(A+\sqrt{4B-A^2})}{3}+Y(-2\sum_{i=1}^{m}\frac{m-i}{3^{i+1}}+\sum_{i=1}^{m}3^{-i}),q_i+m(\frac{Y+A+\sqrt{4B-A^2}}{3})\biggr] $\\
\noindent where $A=(k-2),B=(k+1)^2,Y=(2k+2).$
\end{mydef-corollary}
\begin{proof}
Let $\mathbf{1}_{n}$ be the all $1$ vector of dimension $n\times 1$. Let $Y_1,\dots,Y_n$ be the eigenvectors corresponding to seed graph's signless laplacian eigenvalues $q_1,\dots,q_n$. Before presenting the eigenvectors corresponding to the part (a) of the corollary, let $P^{(i)}_{l_i}$ be the eigenvector corresponding to $G^{(i)}$ (where $i\in[1,m]$) and $l_i\in[1,3n(n+1)^{i-1}]$ such that $\Lambda=\frac{q_i-k+1}{q_i^2-q_i(k+2)+(k+1)}$ be the signless laplacian eigenvalue (defined by above expression) corresponding to $P^{(i)}_{l_i}$. The eigenvectors are as follows
$$
P^{(m)}_{l_m}=\begin{pmatrix}
			P^{(m-1)}_{l_{m-1}}\\
			\Lambda\begin{pmatrix}
										\mathbf{1}_n\otimes P^{(m-1)}_{l_{m-1}}
									 \end{pmatrix}
		\end{pmatrix}
$$
where 
$$
P^{(m-1)}_{l_{m-1}}=\begin{pmatrix}
			P^{(m-2)}_{l_{m-2}}\\
			\Lambda\begin{pmatrix}
											\mathbf{1}_n\otimes P^{(m-2)}_{l_{m-2}}
										 \end{pmatrix}
		  \end{pmatrix},
		  \ldots\,,
P^{(1)}_{l_1}=\begin{pmatrix}
			Y_i\\
			\Lambda\begin{pmatrix}
											\mathbf{1}_n\otimes Y_i
										 \end{pmatrix}
		\end{pmatrix}.
$$

$q_i+1\in\sigma(G^{(1)})$ where for $i=2,\dots,n-1$ with multiplicity of each of them as $k(k+1)^{(m-1)}$.
\end{proof} 
In Fig.~\ref{12}, we plot the distribution of signless Laplacian eigenvalues generated with the seed graph $G^{(0)}=S_3$.

\section{Conclusion}
We proposed a model for generation of complex networks inspired by the phenomena of duplication of genes. We defined Corona graphs by taking corona product of a simple graph, which we call a seed graph, finite number of times. We determined structural properties of the Corona graphs including cumulative degree distribution for any seed graph and cumulative betweenness distribution when the seed graph is a clique. We determined  spectra, Laplacian spectra and signless Laplacian spectra for corona graphs when the seed graph is regular. We also derived the spectra and signless Laplacian spectra of corona graphs when the seed graph is a star graph.\label{sec:Concl}

\bibliography{mybibfile}

\begin{thebibliography}{10}
\expandafter\ifx\csname url\endcsname\relax
  \def\url#1{\texttt{#1}}\fi
\expandafter\ifx\csname urlprefix\endcsname\relax\def\urlprefix{URL }\fi
\expandafter\ifx\csname href\endcsname\relax
  \def\href#1#2{#2} \def\path#1{#1}\fi

\bibitem{leskovec2010kronecker}
J.~Leskovec, D.~Chakrabarti, J.~Kleinberg, C.~Faloutsos, Z.~Ghahramani,
  Kronecker graphs: An approach to modeling networks, The Journal of Machine
  Learning Research 11 (2010) 985--1042.

\bibitem{parsonage2011generalized}
E.~Parsonage, H.~X. Nguyen, R.~Bowden, S.~Knight, N.~Falkner, M.~Roughan,
  Generalized graph products for network design and analysis, in: Network
  Protocols (ICNP), 2011 19th IEEE International Conference on, IEEE, 2011, pp.
  79--88.

\bibitem{frucht1970corona}
R.~Frucht, F.~Harary, On the corona of two graphs, Aequationes Mathematicae
  4~(3) (1970) 322--325.

\bibitem{ispolatov2005duplication}
I.~Ispolatov, P.~Krapivsky, A.~Yuryev, Duplication-divergence model of protein
  interaction network, Physical review E 71~(6) (2005) 061911.

\bibitem{barik2007spectrum}
S.~Barik, S.~Pati, B.~Sarma, The spectrum of the corona of two graphs, SIAM
  Journal on Discrete Mathematics 21~(1) (2007) 47--56.

\bibitem{laurienti2011universal}
P.~J. Laurienti, K.~E. Joyce, Q.~K. Telesford, J.~H. Burdette, S.~Hayasaka,
  Universal fractal scaling of self-organized networks, Physica A: Statistical
  Mechanics and its Applications 390~(20) (2011) 3608--3613.

\bibitem{i2003optimization}
R.~F. i~Cancho, R.~V. Sol{\'e}, Optimization in complex networks, in:
  Statistical mechanics of complex networks, Springer, 2003, pp. 114--126.

\bibitem{zhou2006hierarchical}
C.~Zhou, J.~Kurths, Hierarchical synchronization in complex networks with
  heterogeneous degrees, Chaos: An Interdisciplinary Journal of Nonlinear
  Science 16~(1) (2006) 015104.

\bibitem{wagner2001small}
A.~Wagner, D.~A. Fell, The small world inside large metabolic networks,
  Proceedings of the Royal Society of London. Series B: Biological Sciences
  268~(1478) (2001) 1803--1810.

\bibitem{newman2003structure}
M.~E. Newman, The structure and function of complex networks, SIAM review
  45~(2) (2003) 167--256.

\bibitem{albert2002statistical}
R.~Albert, A.-L. Barab{\'a}si, Statistical mechanics of complex networks,
  Reviews of modern physics 74~(1) (2002) 47.

\bibitem{misiewicz2011fat}
J.~Misiewicz, Fat-tailed distributions: Data, diagnostics, and dependence.

\bibitem{crovella1999estimating}
M.~E. Crovella, M.~S. Taqqu, Estimating the heavy tail index from scaling
  properties, Methodology and computing in applied probability 1~(1) (1999)
  55--79.

\bibitem{rachev2003handbook}
S.~T. Rachev, Handbook of Heavy Tailed Distributions in Finance: Handbooks in
  Finance, Vol.~1, Elsevier, 2003.

\bibitem{jung2002geometric}
S.~Jung, S.~Kim, B.~Kahng, Geometric fractal growth model for scale-free
  networks, Physical Review E 65~(5) (2002) 056101.

\bibitem{zhang2006deterministic}
Z.~Zhang, L.~Rong, C.~Guo, A deterministic small-world network created by edge
  iterations, Physica A: Statistical Mechanics and its Applications 363~(2)
  (2006) 567--572.

\bibitem{boccaletti2006complex}
S.~Boccaletti, V.~Latora, Y.~Moreno, M.~Chavez, D.-U. Hwang, Complex networks:
  Structure and dynamics, Physics reports 424~(4) (2006) 175--308.

\bibitem{holme2002edge}
P.~Holme, Edge overload breakdown in evolving networks, Physical Review E
  66~(3) (2002) 036119.

\bibitem{holme2002vertex}
P.~Holme, B.~J. Kim, Vertex overload breakdown in evolving networks, Physical
  Review E 65~(6) (2002) 066109.

\bibitem{qi2009structural}
Y.~Qi, Z.~Zhang, B.~Ding, S.~Zhou, J.~Guan, Structural and spectral properties
  of a family of deterministic recursive trees: rigorous solutions, Journal of
  Physics A: Mathematical and Theoretical 42~(16) (2009) 165103.

\bibitem{ghim2004packet}
C.-M. Ghim, E.~Oh, K.-I. Goh, B.~Kahng, D.~Kim, Packet transport along the
  shortest pathways in scale-free networks, The European Physical Journal
  B-Condensed Matter and Complex Systems 38~(2) (2004) 193--199.

\bibitem{csardi2006igraph}
G.~Csardi, T.~Nepusz, \href{http://igraph.org}{The igraph software package for
  complex network research}, InterJournal Complex Systems (2006) 1695.
\newline\urlprefix\url{http://igraph.org}

\bibitem{goh2001universal}
K.-I. Goh, B.~Kahng, D.~Kim, Universal behavior of load distribution in
  scale-free networks, Physical Review Letters 87~(27) (2001) 278701.

\bibitem{goh2003goh}
K.-I. Goh, C.-M. Ghim, B.~Kahng, D.~Kim, Goh et al. reply, Physical Review
  Letters 91~(18) (2003) 189804.

\bibitem{vazquez2002large}
A.~V{\'a}zquez, R.~Pastor-Satorras, A.~Vespignani, Large-scale topological and
  dynamical properties of the internet, Physical Review E 65~(6) (2002) 066130.

\bibitem{bapat2010graphs}
R.~B. Bapat, Graphs and matrices, Springer, 2010.

\bibitem{cui2012spectrum}
S.-Y. Cui, G.-X. Tian, The spectrum and the signless laplacian spectrum of
  coronae, Linear Algebra and its Applications 437~(7) (2012) 1692--1703.

\bibitem{cvetkovic2009towards}
D.~Cvetkovi{\'c}, S.~K. Simi{\'c}, Towards a spectral theory of graphs based on
  the signless laplacian, {I}, Publ. Inst. Math.(Beograd) 85~(99) (2009)
  19--33.

\bibitem{haemers2004enumeration}
W.~H. Haemers, E.~Spence, Enumeration of cospectral graphs, European Journal of
  Combinatorics 25~(2) (2004) 199--211.

\bibitem{cvetkovic2010towards}
D.~Cvetkovi{\'c}, S.~K. Simi{\'c}, Towards a spectral theory of graphs based on
  the signless laplacian, {II}, Linear Algebra and its Applications 432~(9)
  (2010) 2257--2272.

\bibitem{cvetkovic2010toward}
D.~Cvetkovi{\'c}, S.~K. Simi{\'c}, Towards a spectral theory of graphs based on
  the signless laplacian, {III}, Applicable Analysis and Discrete Mathematics
  4~(1) (2010) 156--166.

\end{thebibliography}

\end{document}